\documentclass{amsart}

\usepackage{amssymb,amsmath,amsfonts,stmaryrd,amsthm}
\usepackage[sans]{dsfont}
\usepackage[titletoc,title]{appendix}
\usepackage{fouridx,chemarrow}
\usepackage[dvips,cmtip,all]{xy}
\usepackage{setspace,paralist}
\usepackage[linkcolor=black]{hyperref}
\usepackage[letterpaper,left=1.22truein,right=1.22truein,bottom=1.2truein,top=1.4truein]{geometry}
\usepackage{lineno,color,doi,latexsym,relsize}
\usepackage{tikz,tikz-cd,mathabx} 
\usetikzlibrary{matrix,arrows,decorations.pathmorphing}
\modulolinenumbers[5]

\theoremstyle{plain}
\newtheorem{theorem}{Theorem}[section]
\newtheorem{definition}[theorem]{Definition}
\newtheorem{lemma}[theorem]{Lemma}

\newtheorem{proposition}[theorem]{Proposition}
\newtheorem{question}[]{Question}
\newtheorem{problem}[]{Problem}
\newtheorem{corollary}[theorem]{Corollary}
\newtheorem*{nontheorem}{Theorem}{\bf}{\it}

\theoremstyle{remark}
\newtheorem{remark}[theorem]{Remark}
\newtheorem{example}[theorem]{Example}

\numberwithin{equation}{section}

\def\length{\operatorname{length}}
\def\Supph{\operatorname{Supph}}
\def\Supp{\operatorname{Supp}}
\def\Spec{\operatorname{Spec}}
\def\Tor{\operatorname{Tor}}

\def\Perf{\operatorname{\mathbf{Perf}}}
\def\Perfm{\operatorname{\mathbf{Perf}}_{\mathfrak{m}}}
\def\Perfn{\operatorname{\mathbf{Perf}}_{\mathfrak{n}}}

\begin{document}

\title{Entropy in the category of perfect complexes with cohomology of finite length}
\author{Mahdi Majidi-Zolbanin}
\email{mmajidi-zolbanin@lagcc.cuny.edu}
\thanks{The first author received partial funding from a $\mathrm{C}^3\mathrm{IRG}$ (round 10) grant provided by the City University of New York, while working on parts of this paper.}
\thanks{Final version, published in Journal of Pure and Applied Algebra, see:\ \hyperref[https://doi.org/10.1016/j.jpaa.2018.09.008]{https://doi.org/10.1016/j.jpaa.2018.09.008}.}

\author{Nikita Miasnikov}
\email{nikita.miasnikov@oswego.edu}

\address[Mahdi Majidi-Zolbanin]{Department of Mathematics, LaGuardia Community College of the City University of New York, 31-10 Thomson Avenue, Long Island City, NY 11101}
\address[Nikita Miasnikov]{Department of Mathematical Sciences, State University of New York at Oswego, 7060 Route 104, Oswego, NY 13126}


\begin{abstract}
Local and category-theoretical entropies associated with an endomorphism of finite length (i.e., with zero-dimensional closed fiber) of a commutative Noetherian local ring are compared. Local entropy is shown to be less than or equal to category-theoretical entropy. The two entropies are shown to be equal when the ring is regular, and also for the Frobenius endomorphism of a complete local ring of positive characteristic.\par Furthermore, given a flat morphism of Cohen-Macaulay local rings endowed with compatible endomorphisms of finite length, it is shown that local entropy is ``additive''. Finally, over a ring that is a homomorphic image of a regular local ring, a formula for local entropy in terms of an asymptotic partial Euler characteristic is given.
\end{abstract}
\maketitle

\setcounter{section}{-1}
\section{Introduction}
\label{Section:0}
All rings in this paper are assumed to be commutative, Noetherian, and with identity element $1$. Over a commutative ring $R$ we will denote the category of complexes of $R$-modules by $\mathbf{C}(R)$ and the derived category of the category of $R$-modules by $\mathbf{D}(R)$. All homomorphisms of local rings are assumed to be local homomorphisms.\par 
In dynamical systems the complexity of an endomorphism in a given category is usually measured by numerical invariants known as entropy. Often more than one type of entropy may be available to measure the complexity of an endomorphism in a particular category, giving rise to several invariants for the same endomorphism. It is then natural to ask about possible relationships between these invariants. This question has been the focus of many papers. A survey of important results, open problems, and conjectures related to this question, in the category of compact connected Riemannian manifolds can be found in~\cite{LlibSag}. This question is also the main impetus for our work in Section~\ref{Section:2} of this paper, as sketched below:\par Let $(R,\mathfrak{m})$ be a commutative Noetherian local ring. Two types of entropies can be associated to an endomorphism \emph{of finite length} $\phi\colon R\rightarrow R$ (see~Definition~\ref{Definition:1.1}). On one hand there is the local entropy of $\phi$, denoted $h_{\mathrm{loc}}(\phi)$, defined in~\cite[Theorem~1]{MajMiaSzp}. On the other hand, there is a category-theoretical entropy defined in~\cite[Definition~2.1]{Kons} for exact endofunctors of a triangulated category with generator. To associate this type of entropy to $\phi$, we work in $\mathbf{D}(R)$ and note  that the strictly full subcategory of $\mathbf{D}(R)$ formed by perfect complexes with cohomology of finite length, denoted by $\Perfm(R)$, is a triangulated category with generator. Furthermore, the restriction of the total derived inverse image functor $\mathbb{L}\phi^\varstar\colon \mathbf{D}(R)\rightarrow \mathbf{D}(R)$ to $\Perfm(R)$ gives rise to an exact endofunctor of $\Perfm(R)$. This endofunctor has a category-theoretical entropy that is denoted by $h_t(\mathbb{L}\phi^\varstar)$. We should remark that $h_t(\mathbb{L}\phi^\varstar)$ is, by definition, a function of a real variable $t$. In cases that are of particular interest to us, however, $h_t(\mathbb{L}\phi^\varstar)$ turns out to be a constant function. For definitions and details related to above statements, see Section~\ref{Section:1}. \par Section~\ref{Section:2} of this work studies the relationship between the two entropies introduced above. We prove that $h_{\mathrm{loc}}(\phi)\leq h_t(\mathbb{L}\phi^{\varstar})$ for each $t\in\mathbb{R}$, and that equality holds when $R$ is regular, and also when $\phi$ is the Frobenius endomorphism of a \emph{complete} local ring of positive characteristic; see Corollaries~\ref{Corollary:2.2},~\ref{Corollary:2.6}, and Theorem~\ref{Theorem:2.4}.\par
Sections~\ref{Section:3} and~\ref{Section:4} are primarily concerned with further properties of local entropy. Certain invariants of local rings, such as dimension and depth satisfy an ``additivity'' property under flat extensions. That is, given a flat homomorphism $f\colon R\rightarrow S$ of commutative Noetherian local rings, the difference between dimensions (depths) of $S$ and $R$ is equal to the dimension (depth) of the closed fiber of $f$. Our main result in Section~\ref{Section:3}, Theorem~\ref{Theorem:3.3}, is a similar ``additivity'' property for local entropy, under flat extensions of Cohen-Macaulay local rings. To be more precise, given a flat homomorphism $f\colon R\rightarrow S$ of Cohen-Macaulay local rings, and two endomorphisms of finite length $\phi\colon R\rightarrow R$ and $\psi\colon S\rightarrow S$, satisfying $f\circ\phi=\psi\circ f$, we prove $$h_{\mathrm{loc}}(\psi)=h_{\mathrm{loc}}(\phi)+h_{\mathrm{loc}}(\overline{\psi}),$$ where $\overline{\psi}$ is the endomorphism induced by $\psi$ on the closed fiber of $f$.\par
In Section~\ref{Section:4}, Theorem~\ref{Theorem:4.1}, we prove a formula expressing local entropy in terms of an asymptotic partial Euler characteristic, under certain conditions. And in Section~\ref{Section:5} we list a couple of open problems.\par 
\textbf{Acknowledgment.}~We would like to express our gratitude to the referee for his/her detailed comments and suggestions that helped make this work more concise.
\section{Preliminaries}
\label{Section:1}
In this section we recall a number of definitions and basic facts used in this work about local and category-theoretical entropies, as well as perfect complexes. 
\subsection{Local entropy}
\label{Subsection:1.1}
\begin{definition}[\textnormal{\cite[Definition~1]{MajMiaSzp}}]
\label{Definition:1.1}
A local homomorphism \(f\colon R\rightarrow S\) of Noetherian local rings is said to be \emph{of finite length} if its closed fiber is of dimension zero.
\end{definition}
One can quickly see that a local homomorphism \(f\colon(R,\mathfrak{m})\rightarrow(S,\mathfrak{n})\) of Noetherian local rings is of finite length if and only if it satisfies any of the following (equivalent) conditions:
\begin{compactenum}
\item[(a)] \(f(\mathfrak{m})S\) is \(\mathfrak{n}\)-primary;
\item[(b)] If \(\mathfrak{p}\) is a prime ideal of \(S\) such that \(f^{-1}(\mathfrak{p})=\mathfrak{m}\), then \(\mathfrak{p}=\mathfrak{n}\);
\item[(c)] If \(\mathfrak{q}\) is any \(\mathfrak{m}\)-primary ideal of \(R\), then \(f(\mathfrak{q})S\) is \(\mathfrak{n}\)-primary.
\end{compactenum}
\begin{definition}[\textnormal{\cite[Definition~5]{MajMiaSzp}}]
\label{Definition:1.2}
A \emph{local algebraic dynamical system} consists of a Noetherian local ring $(R,\mathfrak{m})$ and an endomorphism of finite length $\phi\colon R\rightarrow R$. We will denote this by $(R,\mathfrak{m},\phi)$. A morphism $f\colon(R,\mathfrak{m},\phi)\rightarrow(S,\mathfrak{n},\psi)$ between two local algebraic dynamical systems is a local homomorphism $f\colon R\rightarrow S$ that satisfies the condition $\psi\circ f=f\circ\phi$. 
\end{definition}
\begin{definition}
\label{Definition:1.3}
Let \((R,\mathfrak{m},\phi)\) be a local algebraic dynamical system and let $\mathfrak{q}$ be an $\mathfrak{m}$-primary ideal of $R$. The local entropy of $\phi$ is the real number defined as follows:
$$h_{\mathrm{loc}}(\phi)=\lim_{n\rightarrow\infty}\frac{1}{n}\log\left(\length_R(R/\phi^n(\mathfrak{q})R)\right).$$
\end{definition}
It is shown in~\cite[Theorems~1, 18]{MajMiaSzp} that $h_{\mathrm{loc}}(\phi)$ is well-defined. That is, the limit defining $h_{\mathrm{loc}}(\phi)$ exists, and is independent of the $\mathfrak{m}$-primary ideal used. In fact, local entropy can be calculated using any module of finite length, and is non-negative.
\begin{example}
\label{Example:1.4}
Let $k$ be a field and $R=k\llbracket X_1,\ldots X_d\rrbracket$. Suppose $\xi_1,\ldots,\xi_d$ are positive integers and $\phi\colon R\rightarrow R$ is the endomorphism that maps $X_i\mapsto X_i^{\xi_i}$ for $1\leq i\leq d$. Then $h_{\mathrm{loc}}(\phi)=\sum_{i=1}^d\log(\xi_i)$. Indeed, as a $k$-vector space, $R/\phi^n(\mathfrak{m})R$ has a basis consisting of monomials $X_1^{i_1}\cdots X_d^{i_d}$, where $0\leq i_j<\xi_i^n$. This implies that $\length_R\left(R/\phi^n(\mathfrak{m})R\right)=\prod_{i=1}^d\xi_i^n$, and hence the local entropy of $\phi$ is the stated one. 
\end{example}
\subsection{Category-theoretical entropy}
\label{Subsection:1.2}
Let $\mathbf{T}$ be a triangulated category. Recall that a subcategory of $\mathbf{T}$ is called \emph{thick} if it is triangulated, contains every object isomorphic to any of its objects, and contains all direct summands of its objects (cf.~\cite[Definition~2.1.6, p.~74]{Neem}). An object $G$ of $\mathbf{T}$ is called a (classical) generator if the smallest thick subcategory of $\mathbf{T}$ containing $G$ is equal to $\mathbf{T}$ itself (cf.~\cite[Section~2.1]{BondVan}). To say that  $G$ is a generator of $\mathbf{T}$  is equivalent to saying that for every object $E$ of $\mathbf{T}$ there is an object $E^\prime$ and a tower of distinguished triangles
\begin{equation}
\label{Equation:1.1}
\begin{tikzcd}[column sep=tiny] 
E_0\arrow[rr]& & E_1 \arrow[dl] \arrow[rr] & &E_2 \arrow[dl] \arrow[r] &\cdots\arrow[r] & E_{p-1} \arrow[rr] & & E_p\cong E\bigoplus E^\prime \arrow[dl] \\  
& G[n_1] \arrow[ul, dashed] & & G[n_2] \arrow[ul, dashed] & & \cdots & & G[n_p] \arrow[ul, dashed] & 
\end{tikzcd}
\end{equation}
with $E_0=0$, $p\geq0$ and $n_i\in\mathbb{Z}$.  
\begin{definition}[{\cite[Definition~2.1]{Kons}}]
\label{Definition:1.5}
Let $G$ and $E$ be objects of a triangulated category $\mathbf{T}$. Let $t$ be a real number. To each tower of distinguished triangles of the form~\emph{(\ref{Equation:1.1})} we associate the exponential sum $\sum_{i=1}^p e^{n_it}$. Let $S_t\subset\mathbb{R}$ be the set of all such sums for a given $t$. The complexity of $E$ with respect to $G$ is the function $\delta_t(G,E)\colon\mathbb{R}\rightarrow[0,\infty]$ of $t$, given by $\delta_t(G,E)=\inf S_t$.
\end{definition}
Note that $\delta_t(G,E)=+\infty$ if and only if $E$ does not lie in the thick subcategory generated by $G$. Also if $F$ is an exact functor from $\mathbf{T}$ to another triangulated category, then since exact functors preserve triangles (and hence towers), the following inequality holds:
\begin{equation}
\label{Equation:1.2}
\delta_t\left(F(G),F(E)\right)\leq\delta_t(G,E).
\end{equation}
\begin{definition}[{\cite[Definition~2.5]{Kons}}]
\label{Definition:1.6}
Let $F\colon\mathbf{T}\rightarrow\mathbf{T}$ be a triangulated endofunctor of a triangulated category $\mathbf{T}$ with a generator $G$. The \emph{entropy} of $F$ is the function $h_t(F)\colon\mathbb{R}\rightarrow[-\infty,+\infty)$ of $t$, given by
$$h_t(F)=\lim_{n\rightarrow\infty}\frac{1}{n}\log\delta_t\left(G,F^n(G)\right).$$
\end{definition}
It is shown in~\cite[Lemma~2.5]{Kons} that $h_t(F)$ is well-defined, i.e., the limit defining $h_t(F)$ exists and is independent of the choice of generator $G$. 
\subsection{Perfect complexes with cohomology of finite length}
\label{Subsection:1.3}
In this subsection we have collected a number of definitions and facts about the category of perfect complexes over a commutative ring, and its strictly full subcategory formed by perfect complexes with cohomology of finite length. The main reference for this subsection is~\cite{SGA6}.
\begin{definition} 
\label{Definition:1.7}
Let $R$ be a commutative ring. A \emph{strictly perfect} complex on $R$ is a bounded complex of projective $R$-modules of finite type. 
\end{definition}
The statement that follows is well-known and will be used implicitly in this work: if $P^\mathsmaller{\bullet}$ and $E^\mathsmaller{\bullet}$ are complexes of $R$-modules, with $P^\mathsmaller{\bullet}$ strictly perfect, then the two conditions below are equivalent:
\begin{compactenum}
\item[1)] There exists a quasi-isomorphism $P^\mathsmaller{\bullet}\stackrel{\sim}{\rightarrow} E^\mathsmaller{\bullet}$ in $\mathbf{C}(R)$;
\item[2)] $P^\mathsmaller{\bullet}$ and $E^\mathsmaller{\bullet}$ are isomorphic in $\mathbf{D}(R)$. 
\end{compactenum}
\begin{definition}
\label{Definition:1.8}
Let $R$ be a commutative ring. A complex $E^\mathsmaller{\bullet}$ of $R$-modules is \emph{perfect} if it has a left resolution by a strictly perfect complex, that is, if there exists a strictly perfect complex $P^\mathsmaller{\bullet}$ on $R$ and a quasi-isomorphism $P^\mathsmaller{\bullet}\stackrel{\sim}{\rightarrow} E^\mathsmaller{\bullet}$ in $\mathbf{C}(R)$. Equivalently, $E^\mathsmaller{\bullet}$ is perfect if in $\mathbf{D}(R)$ it is isomorphic to a strictly perfect complex. The category of perfect complexes over $R$, denoted by $\Perf(R)$ hereafter, is the strictly full subcategory of $\mathbf{D}(R)$ formed by perfect complexes.
\end{definition}
It is well-known (cf.~\cite[Expos\'e~I, Propositions~4.10, 4.17]{SGA6}) that $\Perf(R)$ is a thick subcategory of $\mathbf{D}(R)$.\par
Let $f\colon R\rightarrow S$ be a homomorphism of commutative rings. The \emph{inverse image} functor $f^\varstar\colon$\textbf{R-Mod}$\:\rightarrow\:$\textbf{S-Mod} is the functor that sends an $R$-module $E$ to the $S$-module $S\otimes_RE$ (the notation $f^\varstar$ is used in~\cite[II.5.1, p.~82]{Bourb}). This functor gives rise to an exact functor $f^\varstar\colon\mathbf{K}(R)\rightarrow\mathbf{K}(S)$ of the homotopy categories of complexes that sends a complex $E^\mathsmaller{\bullet}$ of $R$-modules to the complex $S\otimes_RE^\mathsmaller{\bullet}$ of $S$-modules. It is well-known (cf.~\cite[Expos\'e~I, Corollaire~4.19.1]{SGA6}) that the total derived inverse image functor $\mathbb{L}f^\varstar\colon\mathbf{D}(R)\rightarrow\mathbf{D}(S)$, by restriction induces a functor  $\Perf(R)\rightarrow\Perf(S)$. We should remark here that the total derived inverse image functor $\mathbb{L}f^\varstar$ was generally only defined as a functor $\mathbf{D}^-(R)\rightarrow\mathbf{D}^-(S)$ in~\cite{SGA6, Hart}. Spaltenstein extended the definition of $\mathbb{L}f^\varstar$ to an exact functor $\mathbf{D}(R)\rightarrow\mathbf{D}(S)$, using $K$-flat complexes (see~\cite[Proposition~6.7]{Spalten}). A complex $E^\mathsmaller{\bullet}$ in $\mathbf{K}(R)$ is $K$-flat if for every acyclic complex $B^\mathsmaller{\bullet}$ in $\mathbf{K}(R)$, the complex $E^\mathsmaller{\bullet}\otimes_R B^\mathsmaller{\bullet}$ is acyclic. Spaltenstein showed that every complex $G^\mathsmaller{\bullet}$ in $\mathbf{K}(R)$ has a left $K$-flat resolution, i.e., there exists a $K$-flat complex $E^\mathsmaller{\bullet}$ and a quasi-isomorphism $E^\mathsmaller{\bullet}\stackrel{\sim}{\rightarrow}G^\mathsmaller{\bullet}$. Moreover, $\mathbb{L}f^\varstar(G^\mathsmaller{\bullet})$ can be computed by applying $f^\varstar$ to any left $K$-flat resolution of $G^\mathsmaller{\bullet}$ (see~\cite[Proposition~5.6]{Spalten}). Since any bounded complex of flat $R$-modules is $K$-flat (see~\cite[\S~4.3, Lemme~1, p.~66]{Bourb2}), if an object $G^\mathsmaller{\bullet}$ of $\mathbf{D}(R)$ has a left resolution $E^\mathsmaller{\bullet}\stackrel{\sim}{\rightarrow}G^\mathsmaller{\bullet}$ by a bounded flat complex $E^\mathsmaller{\bullet}$ (perfect complexes, for instance), then $\mathbb{L}f^\varstar(G^\mathsmaller{\bullet})$ can be represented by $f^\varstar(E^\mathsmaller{\bullet})$. 
\begin{definition}[{cf.~\cite[p.~157]{Foxby}}]
\label{Definition:1.9}
Let $R$ be a commutative ring and $E^\mathsmaller{\bullet}$ a complex of $R$-modules. The \emph{cohomological support} of $E^\mathsmaller{\bullet}$ is the subspace $\Supph(E^\mathsmaller{\bullet})\subseteq\Spec R$ of those prime ideals $\mathfrak{p}\in\Spec R$ at which the complex $E_{\mathfrak{p}}^\mathsmaller{\bullet}$ of $R_{\mathfrak{p}}$-modules is not acyclic. Equivalently, $\Supph(E^\mathsmaller{\bullet})=\bigcup_{n\in\mathbb{Z}}\Supp H^n(E^\mathsmaller{\bullet})$, as $H^i(E^\mathsmaller{\bullet}\otimes_RR_{\mathfrak{p}})\cong H^i(E^\mathsmaller{\bullet})\otimes_RR_{\mathfrak{p}}$, for all $i$ \emph{(see~\cite[\S~4.2, Corollaire~2, p.~66]{Bourb2})}.
\end{definition}
Over a Noetherian local ring $(R,\mathfrak{m})$ we will denote by $\Perfm(R)$ the strictly full subcategory of $\Perf(R)$ formed by perfect complexes $E^\mathsmaller{\bullet}$ with $\Supph E^\mathsmaller{\bullet}\subseteq\{\mathfrak{m}\}$. That is, a perfect complex $E^\mathsmaller{\bullet}$ is an object of $\Perfm(R)$ if and only if $H^n(E^\mathsmaller{\bullet})$ is an $R$-module of finite length for every $n\in\mathbb{Z}$. One can quickly verify (and it is well-known) that $\Perfm(R)$ is a thick subcategory of $\Perf(R)$. Furthermore, every nonzero object  in $\Perfm(R)$ is a generator, in the sense defined in Section~\ref{Subsection:1.2}. This follows from the following result, first proved in~\cite[Proof of Theorem~11]{Hop} (see also~\cite[Lemma~1.2]{Neem0}): ``Let $R$ be a commutative Noetherian ring and let $E^\mathsmaller{\bullet}, G^\mathsmaller{\bullet}\in\Perf(R)$ be two perfect complexes. If $\Supph(E^\mathsmaller{\bullet})\subseteq\Supph(G^\mathsmaller{\bullet})$, then $E^\mathsmaller{\bullet}$ is in the smallest thick subcategory of $\Perf(R)$ containing $G^\mathsmaller{\bullet}$.'' Thus, if $G^\mathsmaller{\bullet}\in\Perfm(R)$ is a nonzero object and $\left<G^\mathsmaller{\bullet}\right>$ is the smallest thick subcategory of $\Perf(R)$ containing $G^\mathsmaller{\bullet}$, then $\Perfm(R)\subseteq\left<G^\mathsmaller{\bullet}\right>$. But we also have $\left<G^\mathsmaller{\bullet}\right>\subseteq\Perfm(R)$, as can be checked either directly or using the fact that $\Perfm(R)$ itself is a thick subcategory of $\Perf(R)$.
\begin{proposition}
\label{Proposition:1.10}
Let $f\colon(R,\mathfrak{m})\rightarrow(S,\mathfrak{n})$ be a homomorphism of finite length of Noetherian local rings. Then $\mathbb{L}f^\varstar\colon\Perf(R)\rightarrow\Perf(S)$, by restriction induces an exact functor $\mathbb{L}f^\varstar\colon\Perfm(R)\rightarrow\Perfn(S)$.
\end{proposition}
\begin{proof}
Let $E^\mathsmaller{\bullet}$ be an object of $\Perfm(R)$. Then $\mathbb{L}f^\varstar(E^\mathsmaller{\bullet})$ is an object of $\Perf(S)$. We need to show that $\Supph(\mathbb{L}f^\varstar(E^\mathsmaller{\bullet}))\subseteq\{\mathfrak{n}\}$. Let $\mathfrak{q}\in\Spec(S)$ be a non maximal prime ideal and let $\mathfrak{p}=f^{-1}(\mathfrak{q})$. Then $\mathfrak{p}\neq\mathfrak{m}$, as $f$ is a homomorphism of finite length. We have
$$(\mathbb{L}f^\varstar(E^\mathsmaller{\bullet}))_{\mathfrak{q}}=(S\otimes_R^{\mathbb{L}}E^\mathsmaller{\bullet})_{\mathfrak{q}}\cong S_{\mathfrak{q}}\otimes_{R_{\mathfrak{p}}}^{\mathbb{L}}E_{\mathfrak{p}}^\mathsmaller{\bullet}.$$ As the complex $E_{\mathfrak{p}}^\mathsmaller{\bullet}$ is an acyclic object of $\Perf(R_\mathfrak{p})$, it follows quickly from the K\"unneth Formula (see~\cite[\S~4.7, Corollaire~4, p.~79]{Bourb2}) that the complex $S_{\mathfrak{q}}\otimes_{R_{\mathfrak{p}}}^{\mathbb{L}}E_{\mathfrak{p}}^\mathsmaller{\bullet}$ is also acyclic. This shows that $\Supph(f^\varstar(P^\mathsmaller{\bullet}))\subseteq\{\mathfrak{n}\}$, as wanted.
\end{proof}
In the rest of this paper we will refer to the functor $\mathbb{L}f^\varstar\colon\Perfm(R)\rightarrow\Perfn(S)$ described in Proposition~\ref{Proposition:1.10} as \emph{the exact functor induced by} $f$.
\section{Relationships between local and category-theoretical entropies}
\label{Section:2}
\begin{lemma}
\label{Lemma:2.1}
Let $(R,\mathfrak{m})$ be a Noetherian local ring. Let $G^\mathsmaller{\bullet}\in\Perfm(R)$ be a generator.  By the definition of $\Perf(R)$ there exists a smallest non negative integer $N$ such that $H^j(G^\mathsmaller{\bullet})=0$, for $\mid j\mid>N$. Let $$B:=\max\{\length_R\left(H^j(G^\mathsmaller{\bullet})\right)\mid -N\leq j\leq N\}.$$ Then for any object $E^\mathsmaller{\bullet}$ in $\Perfm(R)$, any integer $\ell$, and any real number $t$:
$$\length_R\left(H^\ell(E^\mathsmaller{\bullet})\right)\leq Be^{\ell t}e^{N\mid t\mid}\cdot\delta_t(G^\mathsmaller{\bullet},E^\mathsmaller{\bullet}).$$
\end{lemma}
\begin{proof}
As $H^0(-)$ is a \emph{cohomological} functor (see, e.g.,~\cite[Definition~1.5.2, p.~39]{Kashpira}), it quickly follows that for any distinguished triangle $X^\mathsmaller{\bullet}\rightarrow Y^\mathsmaller{\bullet}\rightarrow Z^\mathsmaller{\bullet}\rightarrow X^\mathsmaller{\bullet}[1]$ in $\Perfm(R)$ and any integer $\ell$: 
$$\length_R\left(H^\ell(Y^\mathsmaller{\bullet})\right)\leq\length_R\left(H^\ell(X^\mathsmaller{\bullet})\right)+\length_R\left(H^\ell(Z^\mathsmaller{\bullet})\right).$$ 
Using this inequality one can immediately check that in each tower of distinguished triangles of $\Perfm(R)$ for $E^\mathsmaller{\bullet}$, of the form displayed in (\ref{Equation:1.1}), for any integer $\ell$:
\begin{eqnarray*}
\length_R\left(H^\ell(E^\mathsmaller{\bullet})\right)&\leq&\sum_{i=1}^p\length_R\left(H^\ell(G^\mathsmaller{\bullet}[n_i])\right)\\
&=&\sum_{i=1}^p\length_R\left(H^{\ell+n_i}(G^\mathsmaller{\bullet})\right).
\end{eqnarray*}
Let $S_\ell:=\{i\in\mathbb{N}\mid -N\leq \ell+n_i\leq N\}$. Then
\begin{equation}
\label{Equation:2.1}
\length_R\left(H^\ell(E^\mathsmaller{\bullet})\right)\leq\sum_{i=1}^p\length_R\left(H^{\ell+n_i}(G^\mathsmaller{\bullet})\right)\leq B\mid S_\ell\mid.
\end{equation}
Next, noting that $e^x\geq e^{-\mid x\mid}$ for any real number $x$, we have
\begin{eqnarray*}
\sum_{i=1}^pe^{(\ell+n_i)t}&\geq&\sum_{i=1}^pe^{-\mid (\ell+n_i)t\mid}\\
&\geq&\sum_{i\in S_\ell}e^{-\mid (\ell+n_i)t\mid}\\ 
&\geq& e^{-N\mid t\mid}\mid S_\ell\mid.
\end{eqnarray*}
Combining this inequality with~(\ref{Equation:2.1}), we obtain 
$$\length\left(H^\ell(E^\mathsmaller{\bullet})\right)\leq Be^{\ell t}e^{N\mid t\mid}\cdot\sum_{i=1}^pe^{n_it}.$$
As $\delta_t(G^\mathsmaller{\bullet},E^\mathsmaller{\bullet})=\inf\left\{\sum_{i=1}^pe^{n_it}\mid\mathrm{the}\ {n_i}^,\mathrm{s}\ \mathrm{appear\ in\ a\ tower\ of\ the\ form~(\ref{Equation:1.1})} \right\}$, the conclusion follows.
\end{proof}
\begin{corollary}
\label{Corollary:2.2}
Let $(R,\mathfrak{m},\phi)$ be a local algebraic dynamical system (see Definition~\ref{Definition:1.2}) and let $\mathbb{L}\phi^{\varstar}\colon \Perfm(R)\rightarrow\Perfm(R)$ be the exact functor induced by $\phi$. Then the inequality $h_{\mathrm{loc}}(\phi)\leq h_t(\mathbb{L}\phi^{\varstar})$ holds for any real number $t$.
\end{corollary}
\begin{proof}
Let $\{x_1,\ldots,x_d\}$ be a system of parameters of $R$ and let $\mathfrak{q}$ be the ideal of $R$ that they generate. Let $G^\mathsmaller{\bullet}(\underline{\mathbf{x}})$ be the Koszul complex over $R$ constructed from $x_1,\ldots,x_d$. (The nonzero modules in this complex are situated in degrees $-d$ to $0$.) As $G^\mathsmaller{\bullet}(\underline{\mathbf{x}})$ is a bounded complex of free modules, for any positive integer $n$ the complex $\mathbb{L}{\phi^n}^{\varstar}(G^\mathsmaller{\bullet}(\underline{\mathbf{x}}))$ can be represented by $\phi^{n\varstar}\left(G^\mathsmaller{\bullet}(\underline{\mathbf{x}})\right)$, which is the Koszul complex $G^\mathsmaller{\bullet}\left(\phi^n(\underline{\mathbf{x}})\right)$ over $R$, constructed from $\phi^n(x_1),\ldots,\phi^n(x_d)$. Thus,  $H^0(\mathbb{L}{\phi^n}^{\varstar}(G^\mathsmaller{\bullet}(\underline{\mathbf{x}})))=R/\phi^n(\mathfrak{q})R$. We take $G^\mathsmaller{\bullet}(\underline{\mathbf{x}})$ as a generator for the triangulated category $\Perfm(R)$ and apply Lemma~\ref{Lemma:2.1} with $\ell=0$ and $\mathbb{L}{\phi^n}^{\varstar}(G^\mathsmaller{\bullet}(\underline{\mathbf{x}}))$ as $E^\mathsmaller{\bullet}$, to obtain
\begin{equation}
\label{Equation:2.2}
\length_R\left(R/\phi^n(\mathfrak{q})R\right)\leq Be^{N\mid t\mid}\cdot\delta_t\left(G^\mathsmaller{\bullet}(\underline{\mathbf{x}}),\mathbb{L}{\phi^n}^{\varstar}(G^\mathsmaller{\bullet}(\underline{\mathbf{x}}))\right),
\end{equation}
where $B$ and $N$ are constants defined in that lemma. Now the desired inequality $h_{\mathrm{loc}}(\phi)\leq h_t(\mathbb{L}\phi^{\varstar})$ follows by taking the logarithm, dividing by $n$, and passing to the limit as $n\rightarrow\infty$ on both sides of (\ref{Equation:2.2}).
\end{proof}
\begin{remark}
\label{Remark:2.3}
In a Cohen-Macaulay Noetherian local ring of dimension $d$, a sequence of $d$ elements in the maximal ideal form a system of parameters if and only if they form a (maximal) regular sequence. For a proof of this fact see~\cite[Theorem~17.4]{Matsumura2}. We will use this fact a few times in this paper, for instance in proofs of  Theorems~\ref{Theorem:2.4} and~\ref{Theorem:3.3}. 
\end{remark}
\begin{theorem}
\label{Theorem:2.4}
Let $(R,\mathfrak{m},\phi)$ be a local algebraic dynamical system (see Definition~\ref{Definition:1.2}). Assume that $R$ is regular of dimension $d$, and let $\mathbb{L}\phi^{\varstar}\colon\Perfm(R)\rightarrow\Perfm(R)$ be the exact functor induced by $f$.  Then $h_t(\mathbb{L}\phi^{\varstar})$ is constant and equal to $h_{\mathrm{loc}}(\phi)$.
\end{theorem}
\begin{proof}
As $R$ is regular, every $R$-module of finite type has finite projective dimension and therefore, considered as a complex concentrated in degree zero, is an object of $\Perf(R)$. Let $k=R/\mathfrak{m}$ be the residue field of $R$. We make two claims:
\begin{compactenum}
\item[Claim 1:] if $E$ is an $R$-module of finite type, then $\delta_t(k,E)\leq\length_R(E)$;
\item[Claim 2:] if $n\geq0$, then $\mathbb{L}\phi^{n\varstar}(k)$ can be represented by $R/\phi^n(\mathfrak{m})R$.
\end{compactenum}
Let us first prove the theorem assuming these claims: we take $k$ as generator for the triangulated category $\Perfm(R)$. Using Claims 1 and 2 above, for any integer $n\geq0$ we can write
\begin{eqnarray*}
\delta_t\left(k,\mathbb{L}\phi^{n\varstar}(k)\right)&=&\delta_t(k,R/\phi^n(\mathfrak{m})R)\\
& \leq & \length_R(R/\phi^n(\mathfrak{m})R).
\end{eqnarray*}
Taking the logarithm, dividing by $n$, and passing to the limit as $n\rightarrow\infty$ in the previous inequality, we get $h_t(\mathbb{L}\phi^{\varstar})\leq h_{\mathrm{loc}}(\phi)$. On the other hand, Corollary~\ref{Corollary:2.2} gives us the reverse inequality $h_{\mathrm{loc}}(\phi)\leq h_t(\mathbb{L}\phi^{\varstar})$. Thus, $h_t(\mathbb{L}\phi^{\varstar}) = h_{\mathrm{loc}}(\phi)$. We now prove the claims:\par Proof of Claim 1: If $\length_R(E)=\infty$ then the claim holds trivially. Assume $\length_R(E)<\infty$. We will use induction on $\length_R(E)$. The claim clearly holds if $\length_R(E)=1$, as $0\rightarrow k\rightarrow k\rightarrow0$ is a distinguished triangle in $\Perfm(R)$, showing that $\delta_t(k,k)\leq1$. Suppose now that $\length_R(E)>1$. Then there is an exact sequence of $R$-modules $$0\rightarrow E_1\rightarrow E\rightarrow k\rightarrow 0$$ with $\length_R(E_1)=\length_R(E)-1$. This exact sequence gives rise to a distinguished triangle $E_1\rightarrow E\rightarrow k\rightarrow E_1[1]$ in $\Perfm(R)$ (cf.~\cite[Proposition~1.7.5, p.~46]{Kashpira}). Attaching this distinguished triangle (or its direct sum with a distinguished triangle of the form $E^\prime\rightarrow E^\prime\rightarrow 0\rightarrow E^\prime[1]$, if necessary) to the right end of any tower of distinguished triangles for $E_1$ of the form displayed in~(\ref{Equation:1.1}), will get us a tower of distinguished triangles for $E$, from which it is clear that $\delta_t(k,E)\leq\delta_t(k,E_1)+1$. The claim now follows from the induction hypothesis.\par
Proof of Claim 2: let $\{x_1,\ldots,x_d\}$ be a regular system of parameters of $R$, that is, a set of $d$ elements that generate the maximal ideal $\mathfrak{m}$. Let $G^\mathsmaller{\bullet}(\underline{\mathbf{x}})$ be the Koszul complex over $R$ constructed from $x_1,\ldots,x_d$ (the nonzero modules in this complex are situated in degrees $-d$ to $0$). Considering $k$ as a complex concentrated in degree zero, there is a quasi-isomorphism $G^\mathsmaller{\bullet}(\underline{\mathbf{x}})\stackrel{\sim}{\rightarrow}k$. Hence, for any positive integer $n$ the complex $\mathbb{L}\phi^{n\varstar}(k)$ can be represented by $\phi^{n\varstar}(G^\mathsmaller{\bullet}(\underline{\mathbf{x}}))$, which is the Koszul complex $G^\mathsmaller{\bullet}\left(\phi^n(\underline{\mathbf{x}})\right)$ over $R$, constructed from $\phi^n(x_1),\ldots,\phi^n(x_d)$. As $\phi$ (and hence $\phi^n$) is of finite length, the ideal generated by $\phi^n(x_1),\ldots,\phi^n(x_d)$ is $\mathfrak{m}$-primary, i.e., $\{\phi^n(x_1),\ldots,\phi^n(x_d)\}$ is a system of parameters of $R$. By Remark~\ref{Remark:2.3} then, $\phi^n(x_1),\ldots,\phi^n(x_d)$ is a regular sequence. Thus, $H^i(G^\mathsmaller{\bullet}\left(\phi^n(\underline{\mathbf{x}})\right))=0$ for $i\neq0$, and $H^0(G^\mathsmaller{\bullet}\left(\phi^n(\underline{\mathbf{x}})\right))= R/\phi^n(\mathfrak{m})R$. Hence, considering $R/\phi^n(\mathfrak{m})R$ as a complex concentrated in degree zero, there is a quasi-isomorphism $$G^\mathsmaller{\bullet}\left(\phi^n(\underline{\mathbf{x}})\right)\stackrel{\sim}{\rightarrow}R/\phi^n(\mathfrak{m})R.$$ That is, $\mathbb{L}\phi^{n\varstar}(k)$ can also be represented by $R/\phi^n(\mathfrak{m})R$, as claimed.
\end{proof} 
\begin{proposition}
\label{Proposition:2.5}
Suppose $f\colon(R,\mathfrak{m},\phi)\rightarrow(S,\mathfrak{n},\psi)$ is a morphism of local algebraic dynamical systems, with $f\colon R\rightarrow S$ of finite length. Let $\mathbb{L}\phi^{\varstar}\colon\Perfm(R)\rightarrow\Perfm(R)$ and $\mathbb{L}\psi^{\varstar}\colon\Perfn(S)\rightarrow \Perfn(S)$ be the exact functors induced by $\phi$ and $\psi$, respectively. Then:
\begin{compactenum}
\item[\emph{a)}] $h_t(\mathbb{L}\psi^{\varstar})\leq h_t(\mathbb{L}\phi^{\varstar})$.
\item[\emph{b)}] If in addition $R$ is regular and $h_{\mathrm{loc}}(\phi)=h_{\mathrm{loc}}(\psi)$, then $h_t(\mathbb{L}\psi^{\varstar})$ is constant and equal to $h_{\mathrm{loc}}(\psi)$.
\end{compactenum}
\end{proposition}
\begin{proof}
a) Let $\{x_1,\ldots,x_d\}$ be a system of parameters of $R$, where $d=\dim R$, and let $y_i=f(x_i)$ for $1\leq i\leq d$. Let $G_R^\mathsmaller{\bullet}(\underline{\mathbf{x}})$ and $G_S^\mathsmaller{\bullet}(\underline{\mathbf{y}})$ be the Koszul complexes over $R$ and $S$, respectively, constructed from $x_1,\ldots,x_d$ and $y_1,\ldots,y_d$. We take $G_R^\mathsmaller{\bullet}(\underline{\mathbf{x}})$ and $G_S^\mathsmaller{\bullet}(\underline{\mathbf{y}})$ as generators of the triangulated categories $\Perfm(R)$ and $\Perfn(S)$, respectively. Let $\mathbb{L}f^{\varstar}\colon\Perfm(R)\rightarrow \Perfn(S)$ be the exact functor induced by $f$. As $f^{\varstar}(G_R^\mathsmaller{\bullet}(\underline{\mathbf{x}}))=G_S^\mathsmaller{\bullet}(\underline{\mathbf{y}})$ and $\mathbb{L}f^{\varstar}(G_R^\mathsmaller{\bullet}(\underline{\mathbf{x}}))$ can be represented by the complex $f^{\varstar}(G_R^\mathsmaller{\bullet}(\underline{\mathbf{x}}))$, we can write 
\begin{equation}
\label{Equation:2.3}
\mathbb{L}f^{\varstar}(G_R^\mathsmaller{\bullet}(\underline{\mathbf{x}}))=G_S^\mathsmaller{\bullet}(\underline{\mathbf{y}}).
\end{equation} 
The condition $f\circ\phi=\psi\circ f$ satisfied by $f$ for being a morphism of local algebraic dynamical systems gives us 
\begin{equation}
\label{Equation:2.4}
\mathbb{L}f^{\varstar}\circ\mathbb{L}\phi^{\varstar}=\mathbb{L}\psi^{\varstar}\circ\mathbb{L}f^{\varstar}.
\end{equation} 
Now for any integer $n\geq1$ and any real number $t$, using equalities~(\ref{Equation:2.3}) and~(\ref{Equation:2.4}) we can write:
\begin{eqnarray*}
\delta_t\big(G_S^\mathsmaller{\bullet}(\underline{\mathbf{y}}),\mathbb{L}{\psi^n}^\varstar(G_S^\mathsmaller{\bullet}(\underline{\mathbf{y}}))\big)&=&\delta_t\left(\mathbb{L}f^{\varstar}\left(G_R^\mathsmaller{\bullet}(\underline{\mathbf{x}})\right),\mathbb{L}{\psi^n}^{\varstar}\left(\mathbb{L}f^{\varstar}\left(G_R^\mathsmaller{\bullet}(\underline{\mathbf{x}})\right)\right)\right)\\
&=&\delta_t\left(\mathbb{L}f^{\varstar}\left(G_R^\mathsmaller{\bullet}(\underline{\mathbf{x}})\right),\mathbb{L}f^{\varstar}\left(\mathbb{L}{\phi^n}^{\varstar}\left(G_R^\mathsmaller{\bullet}(\underline{\mathbf{x}})\right)\right)\right)\\
&\leq&\delta_t\left(G_R^\mathsmaller{\bullet}(\underline{\mathbf{x}}),\mathbb{L}{\phi^n}^{\varstar}\left(G_R^\mathsmaller{\bullet}(\underline{\mathbf{x}})\right)\right),
\end{eqnarray*}
where the last inequality holds by~(\ref{Equation:1.2}). By taking the logarithm, dividing by $n$, and passing to the limit as $n\rightarrow\infty$ we obtain $h_t(\mathbb{L}\psi^{\varstar})\leq h_t(\mathbb{L}\phi^{\varstar})$.\par b) Combining part a) with the result of Corollary~\ref{Corollary:2.2} we obtain:
$$h_{\mathrm{loc}}(\psi)\leq h_t(\mathbb{L}\psi^{\varstar})\leq h_t(\mathbb{L}\phi^{\varstar}).$$
If $R$ is regular, then $h_t(\mathbb{L}\phi^{\varstar})=h_{\mathrm{loc}}(\phi)$ by Theorem~\ref{Theorem:2.4}. Since $h_{\mathrm{loc}}(\phi)=h_{\mathrm{loc}}(\psi)$ by assumption, we conclude that $h_t(\mathbb{L}\psi^{\varstar})$ is constant and equal to $h_{\mathrm{loc}}(\psi)$.
\end{proof}
\begin{corollary}
\label{Corollary:2.6}
Let $(S,\mathfrak{n})$ be an arbitrary complete Noetherian local ring of positive characteristic $p$ and dimension $d$, let $f_S\colon S\rightarrow S$ be the Frobenius endomorphism of $S$, and let $\mathbb{L}f_S^{\varstar}\colon\Perfn(S)\rightarrow\Perfn(S)$ be the exact functor induced by $f_S$. Then $h_t(\mathbb{L}f_S^{\varstar})$ is constant and equal to $d\cdot\log(p)$.
\end{corollary}
\begin{proof}
Let $\{x_1,\ldots,x_d\}$ be a system of parameters of $S$, and $k$ the residue field of $S$. Recall that $S$ is a module-finite extension of the regular ring $R:=k\llbracket X_1,\ldots,X_d\rrbracket$ via the injective ring homomorphism $\eta\colon R\rightarrow S$ that maps $X_i$ onto $x_i$, for $1\leq i\leq d$ (cf.~\cite[Theorem~29.4, p.~225]{Matsumura2}). Let $f_R$ be the Frobenius endomorphism of $R$. By~\cite[Theorem~1]{MajMiaSzp} the local entropy of the Frobenius endomorphism of a Noetherian local ring of characteristic $p>0$ and of dimension $d$ is equal to $d\cdot\log(p)$. Thus, $h_{\mathrm{loc}}(f_R)=h_{\mathrm{loc}}(f_S)=d\cdot\log p$. Since $\eta\circ f_R=f_S\circ\eta$, the result follows from Proposition~\ref{Proposition:2.5}.
\end{proof}
\section{Additivity of local entropy under flat extensions} 
\label{Section:3}
Certain invariants of local rings, such as dimension and depth, are ``additive'' under flat extensions. That is, if $f\colon(R,\mathfrak{m})\rightarrow S$ is a flat homomorphism of commutative Noetherian local rings, then
\begin{equation}
\label{Equation:3.1}
\dim S=\dim R+\dim S/f(\mathfrak{m})S,
\end{equation} and the same equation holds replacing dimension with depth. Craig~Huneke asked us whether local entropy satisfies a similar ``additivity'' property under flat extensions. To be more precise, let \(f\colon(R,\mathfrak{m},\phi)\rightarrow(S,\mathfrak{n},\psi)\) be a morphism of local algebraic dynamical systems. Then by definition of such morphisms, the relation $\psi\circ f=f\circ\phi$ holds, from which it quickly follows that the ideal $f(\mathfrak{m})S$ is $\psi$-stable, that is, $$\psi\left(f(\mathfrak{m})S\right)\subseteq f(\mathfrak{m})S.$$ Thus, $\psi$ induces an endomorphism of finite length $\overline{\psi}\colon S/f(\mathfrak{m})S\rightarrow S/f(\mathfrak{m})S$ on the closed fiber of $f$. Under these settings, Huneke's question can be formulated as follows:
\begin{question}
\label{Question:1}
If $f$ is flat, does it hold that $h_{\mathrm{loc}}(\psi)=h_{\mathrm{loc}}(\phi)+h_{\mathrm{loc}}(\overline{\psi})?$
\end{question}
If $\dim R=\dim S$, then Question~\ref{Question:1} has an affirmative answer. This is proved in~\cite[Corollary~16 and Proposition~20]{MajMiaSzp}. Question~\ref{Question:1} has also an affirmative answer when $\phi$ and $\psi$, respectively, are the Frobenius endomorphisms of two local rings $R$ and $S$ of characteristic $p>0$. Indeed, as the local entropy of the Frobenius endomorphism of a local ring of characteristic $p>0$ and of dimension $d$ is equal to $d\cdot\log p$ (see~\cite[Theorem~1]{MajMiaSzp}), in this case the equality in Question~\ref{Question:1} quickly reduces to (\ref{Equation:3.1}), which holds, since $f$ is flat (see, e.g.,~\cite[Theorem~15.1]{Matsumura2}). \par Our main goal in this section is to give an affirmative answer to Question~\ref{Question:1}, in Theorem~\ref{Theorem:3.3}, in the special case when $S$ is Cohen-Macaulay. The question remains open in the general non-Cohen-Macaulay case.\par
We will use the following Flatness Criterion in the proof of Theorem~\ref{Theorem:3.3}, as well as in Example~\ref{Example:3.4}. See~\cite[Corollary to Theorem~22.5]{Matsumura2} for a proof of this criterion.
\begin{nontheorem}[Flatness Criterion]
\textit{Let \(f\colon(R,\mathfrak{m})\rightarrow(S,\mathfrak{n})\) be a local homomorphism of Noetherian local rings and let $M$ be an $S$-module of finite type. For $y_1,\ldots,y_n\in\mathfrak{n}$ write $\overline{y}_i$ for the images of $y_i$ in $S/f(\mathfrak{m})S$. Then the following conditions are equivalent:}
\begin{compactenum}
\item[a)] \textit{$y_1,\ldots,y_n$ is an $M$-regular sequence and $M/\sum_1^ny_iM$ is flat over $R$;}
\item[b)] \textit{$\overline{y}_1,\ldots,\overline{y}_n$ is an $(M/f(\mathfrak{m})M)$-regular sequence and $M$ is flat over $R$.}
\end{compactenum}
\end{nontheorem}
We will also need the following elementary statement:
\begin{proposition}
\label{Proposition:3.1}
Let \(f:(R,\mathfrak{m})\rightarrow S\) be a local homomorphism of finite length of Noetherian local rings. Let \(M\) be an \(R\)-module of finite length. Then
\begin{compactenum}
\item[a)] \(M\otimes_RS\) is of finite length as an \(S\)-module.
\item[b)] \(\length_S(M\otimes_R S)\leq\length_R(M)\cdot\length_S(S/f(\mathfrak{m})S)\).
\item[c)] If \(f\) is \textit{flat}, then \(\length_S(M\otimes_R S)=\length_R(M)\cdot\length_S(S/f(\mathfrak{m})S)\).
\end{compactenum}
\end{proposition}
\begin{proof}
By induction on $\length_R(M)$.
\end{proof}
We begin with showing that a morphism of local algebraic dynamical systems gives rise to an inequality between local entropies:
\begin{proposition}
\label{Proposition:3.2}
Suppose \(f\colon(R,\mathfrak{m},\phi)\rightarrow(S,\mathfrak{n},\psi)\) is a morphism of local algebraic dynamical systems and let $\overline{\psi}\colon S/f(\mathfrak{m})S\rightarrow S/f(\mathfrak{m})S$ be the endomorphism induced by $\psi$ (see the paragraph before Question~\ref{Question:1}). Then the following inequality holds:
$$h_{\mathrm{loc}}(\psi)\leq h_{\mathrm{loc}}(\phi)+h_{\mathrm{loc}}(\overline{\psi}).$$
\end{proposition}
\begin{proof}
The composition of maps $R\stackrel{f}{\rightarrow}S\rightarrow S/\psi^n(\mathfrak{n})S$ gives a local homomorphism of finite length $R\rightarrow S/\psi^n(\mathfrak{n})S$ for each integer $n\geq0$. Applying Proposition~\ref{Proposition:3.1}, we can write:
\begin{eqnarray*}
\length_S(S/\psi^n(\mathfrak{n})S)&=&\length_S\left((R/\phi^n(\mathfrak{m})R)\otimes_R (S/\psi^n(\mathfrak{n})S\right)\\ &\leq&\length_R(R/\phi^n(\mathfrak{m})R)\cdot\length_S(S/(f(\mathfrak{m})S+\psi^n(\mathfrak{n})S)).
\end{eqnarray*}
We obtain the desired inequality by applying logarithm, dividing by $n$ and taking limits as $n\rightarrow\infty$.
\end{proof}

We now give an affirmative answer to Question~\ref{Question:1} in the particular case when $S$ is Cohen-Macaulay:
\begin{theorem}
\label{Theorem:3.3}
Suppose \(f\colon(R,\mathfrak{m},\phi)\rightarrow(S,\mathfrak{n},\psi)\) is a flat morphism of local algebraic dynamical systems and let $\overline{\psi}\colon S/f(\mathfrak{m})S\rightarrow S/f(\mathfrak{m})S$ be the endomorphism induced by $\psi$ (see the paragraph before Question~\ref{Question:1}). If $S$ is Cohen-Macaulay, then
\begin{equation}
\label{Equation:3.2}
h_{\mathrm{loc}}(\psi)=h_{\mathrm{loc}}(\phi)+h_{\mathrm{loc}}(\overline{\psi}).
\end{equation}
\end{theorem}
\begin{proof}
As $f$ is flat, the Cohen-Macaulayness of $S$ implies that the rings $R$ and $S/f(\mathfrak{m})S$ are also Cohen-Macaulay (see, e.g.,~\cite[Corollary to Theorem~23.3]{Matsumura2}). Since $S/f(\mathfrak{m})S$ is Cohen-Macaulay, there exists a (non-unique) sequence of elements $y_1,\ldots,y_{d^\prime}\in\mathfrak{n}$ of length $d^\prime=\dim (S/f(\mathfrak{m})S)$, whose images in $S/f(\mathfrak{m})S$ form an $(S/f(\mathfrak{m})S)$-regular sequence. Note that by the Flatness Criterion stated earlier, $y_1,\ldots,y_{d^\prime}$ is an $S$-regular sequence. Let $\mathfrak{q}^\prime\subset S$ be the ideal generated by $y_1,\ldots,y_{d^\prime}$. We claim that for any integer $n\geq0$, the ring $S/\psi^n(\mathfrak{q}^\prime)S$ is flat over $R$ via the composition of maps
\begin{equation}
\label{Equation:3.3}
R\stackrel{f}{\rightarrow}S\rightarrow S/\psi^n(\mathfrak{q}^\prime)S.
\end{equation}
As $R\stackrel{f}{\rightarrow}S$ is flat, the claim will be established by the Flatness Criterion, if we can show that the images of $\psi^n(y_1),\ldots,\psi^n(y_{d^\prime})$ in $S/f(\mathfrak{m})S$ form an $(S/f(\mathfrak{m})S)$-regular sequence. These images coincide with elements $$\overline{\psi}^n(\overline{y}_1),\ldots,\overline{\psi}^n(\overline{y}_{d^\prime}),$$ where $\overline{y}_i$ is the image of $y_i$ in $S/f(\mathfrak{m})S$. That $\overline{\psi}^n(\overline{y}_1),\ldots,\overline{\psi}^n(\overline{y}_{d^\prime})$ is an $(S/f(\mathfrak{m})S)$-regular sequence is an immediate consequence of Remark~\ref{Remark:2.3}, the fact that $\overline{y}_1,\ldots,\overline{y}_{d^\prime}$ is a maximal $(S/f(\mathfrak{m})S)$-regular sequence, and the fact that $\overline{\psi}^n$ is an endomorphism of finite length of $S/f(\mathfrak{m})S$ (hence, the image under $\overline{\psi}^n$ of any system of parameters is again a system of parameters in $S/f(\mathfrak{m})S$).\par
Now let  $x_1,\ldots,x_d\in\mathfrak{m}$ be an $R$-regular sequence of length $d=\dim R$ and let $\mathfrak{q}\subset R$ be the ideal generated by $x_1,\ldots,x_d$. By Remark~\ref{Remark:2.3}, $\mathfrak{q}$ is a parameter ideal of $R$. By the flatness of $S/\mathfrak{q}^\prime$ over $R$ via the composition of maps shown in~(\ref{Equation:3.3}) (taking $n=0$), the images of $f(x_1),\ldots,f(x_d)$ in $S/\mathfrak{q}^\prime$ form an $(S/\mathfrak{q}^\prime)$-regular sequence. This means $y_1,\ldots,y_{d^\prime},f(x_1),\ldots,f(x_d)$ is an $S$-regular sequence. Moreover, since $f$ is flat, $$d+d^\prime=\dim R+\dim (S/f(\mathfrak{m})S)=\dim S$$ (see, e.g.,~\cite[Theorem~15.1]{Matsumura2}). Hence, $\{y_1,\ldots,y_{d^\prime},f(x_1),\ldots,f(x_d)\}$ is a system of parameters in $S$, by Remark~\ref{Remark:2.3}. Let $\mathfrak{Q}\subset S$ be the ideal generated by $$y_1,\ldots,y_{d^\prime},f(x_1),\ldots,f(x_d).$$ We note that for any integer $n\geq0$:
\begin{equation}
\label{Equation:3.4}
\frac{R}{\phi^n(\mathfrak{q})R}\otimes_R\frac{S}{\psi^n(\mathfrak{q}^\prime)S}\cong\frac{S}{f(\phi^n(\mathfrak{q}))S+\psi^n(\mathfrak{q}^\prime)S}\cong\frac{S}{\psi^n(\mathfrak{Q})S},
\end{equation}
where the last isomorphism quickly follows from the fact that $\psi\circ f=f\circ\phi$. Since $S/\psi^n(\mathfrak{q}^\prime)S$ is flat over $R$ and 
$$\dim (S/\psi^n(\mathfrak{q}^\prime)S)=\dim S-d^\prime=\dim S-\dim (S/f(\mathfrak{m})S)=\dim R,$$
the homomorphism $R\rightarrow S/\psi^n(\mathfrak{q}^\prime)S$ obtained by composing the maps given in~(\ref{Equation:3.3}) is in fact, of finite length. Hence, Proposition~\ref{Proposition:3.1}-c) applies and from~(\ref{Equation:3.4}) we obtain
\begin{eqnarray*}
\length_S\left(S/\psi^n(\mathfrak{Q})S\right)&=&\length_S\Big(\frac{R}{\phi^n(\mathfrak{q})R}\otimes_R\frac{S}{\psi^n(\mathfrak{q}^\prime)S}\Big)\\
&=&\length_R\left(R/\phi^n(\mathfrak{q})R\right)\cdot\length_S\left(S/[f(\mathfrak{m})S+\psi^n(\mathfrak{q}^\prime)S]\right).
\end{eqnarray*}
After applying logarithm to both sides, dividing by $n$ and taking limits as $n\rightarrow\infty$, we obtain~(\ref{Equation:3.2}).
\end{proof}
\begin{example}
\label{Example:3.4}
In this example we will apply Theorem~\ref{Theorem:3.3} to calculate local entropy of a specific endomorphism. The local endomorphism of the ring $(\mathbb{Z}/2\mathbb{Z})\llbracket X,Y,W,U\rrbracket$ that maps $X, Y, W$ and $U$ to $X^3+U^3, Y^3, W^5+X^2$ and $XU^2$, respectively, is of finite length, because if $\mathfrak{p}$ is a minimal prime ideal of $(X^3+U^3, Y^3, W^5+X^2,XU^2)$, then as one can quickly see, $\mathfrak{p}=(X,Y,W,U)$. One can also verify quickly that the ideal $(U^6,Y^3+X^2)$ is stable under this endomorphism. Thus, we obtain an induced ring endomorphism of finite length:
$$\psi\colon\frac{(\mathbb{Z}/2\mathbb{Z})\llbracket X,Y,W,U\rrbracket}{(U^6,Y^3+X^2)}\rightarrow\frac{(\mathbb{Z}/2\mathbb{Z})\llbracket X,Y,W,U\rrbracket}{(U^6,Y^3+X^2)}.$$ 
To abbreviate notation we will write $S$ for the ring $(\mathbb{Z}/2\mathbb{Z})\llbracket X,Y,W,U\rrbracket/(U^6,Y^3+X^2)$. Our goal in this example is to calculate $h_{\mathrm{loc}}(\psi)$, the local entropy of $\psi$. We will do this by constructing a flat homomorphism into the ring $S$ and then using Theorem~\ref{Theorem:3.3}. Note that $S$ is Cohen-Macaulay by virtue of being a complete intersection.\par Let $R=(\mathbb{Z}/2\mathbb{Z})\llbracket T\rrbracket$ and let $\phi\colon R\rightarrow R$ be the local endomorphism that maps $T$ to $T^3$. Let $f\colon R\rightarrow S$ be the local homomorphism such that $f(T)=y$, where $y$ is the image of $Y$ in $S$. It is evident that $f\circ\phi=\psi\circ f$. From the Flatness Criterion that was stated earlier, it quickly follows that $f$ is flat. Hence, by Theorem~\ref{Theorem:3.3}
\begin{eqnarray*}
h_{\mathrm{loc}}(\psi)&=&h_{\mathrm{loc}}(\phi)+h_{\mathrm{loc}}(\overline{\psi})\\
&=&\log(3)+h_{\mathrm{loc}}(\overline{\psi}),
\end{eqnarray*}
where as usual $\overline{\psi}$ is the endomorphism induced by $\psi$ on $S/yS$. (That $h_{\mathrm{loc}}(\phi)=\log(3)$ can be calculated quickly, using the definition of local entropy, as seen in Example~\ref{Example:1.4}.) The ring $S/yS$ is isomorphic to $S^\prime:=(\mathbb{Z}/2\mathbb{Z})\llbracket X,W,U\rrbracket/(U^6,X^2)$ and $\overline{\psi}\colon S^\prime\rightarrow S^\prime$ maps $x, w$ and $u$ to $u^3, w^5$ and $xu^2$, respectively, where $x, w$ and $u$ are images of $X, W$ and $U$ in $S^\prime$. In order to calculate $h_{\mathrm{loc}}(\overline{\psi})$, we construct another flat homomorphism, this time into $S^\prime$. Let $R^\prime:=(\mathbb{Z}/2\mathbb{Z})\llbracket Z\rrbracket$ and let $\phi^\prime\colon R^\prime\rightarrow R^\prime$ be the local endomorphism that maps $Z$ to $Z^5$. Let $f^\prime\colon R^\prime\rightarrow S^\prime$ be the local homomorphism such that $f^\prime(Z)=W$. Again it is evident that $f\circ\phi=\psi\circ f$ and the flatness of $f^\prime$ quickly follows from the Flatness Criterion that was stated earlier. By Theorem~\ref{Theorem:3.3}, and using the fact that the local entropy of an endomorphism of a zero-dimensional local ring is zero (\cite[Corollary~16]{MajMiaSzp}), we quickly see that $h_{\mathrm{loc}}(\overline{\psi})=\log(5)$. Hence, $h_{\mathrm{loc}}(\psi)=\log(3)+\log(5)$.
\end{example}
\section{Local entropy as an asymptotic partial Euler characteristic}
\label{Section:4}
When there is a surjective morphism $f\colon(R,\mathfrak{m},\phi)\rightarrow(S,\mathfrak{n},\psi)$ of local algebraic dynamical systems with $R$ regular, then $h_{\mathrm{loc}}(\psi)$, the local entropy of $\psi$, can be expressed as an asymptotic ``partial intersection multiplicity'', as stated in the next theorem.
\begin{theorem}
\label{Theorem:4.1}
Let \(f\colon(R,\mathfrak{m},\phi)\rightarrow(S,\mathfrak{n},\psi)\) be a surjective morphism of local algebraic dynamical systems, that is, $S$ is the homomorphic image of $R$ under $f$. Assume that $\ker f\neq(0)$ and that $R$ is regular of dimension $d$. Then the following equality holds:
\begin{equation}
\label{Equation:4.1}
h_{\mathrm{loc}}(\psi)=\lim_{n\rightarrow\infty}\frac{1}{n}\log\Big(\sum_{i=1}^{d}(-1)^{i-1}\length_R\big(\Tor_i^R\big(R/\phi^n(\mathfrak{m})R,S\big)\big)\Big).
\end{equation}
\end{theorem}
\begin{proof}
The $R$-module $\big(R/\phi^n(\mathfrak{m})R\big)\otimes_RS$ is of finite length and 
$$\dim (R/\phi^n(\mathfrak{m})R) + \dim S = \dim S < \dim R.$$ By the vanishing part of Serre's intersection multiplicity~\cite[Theorem~1, p.~106]{Serre} proven for arbitrary regular local rings in~\cite{Roberts1},~\cite{Roberts2} and independently in~\cite{GillSoul1},~\cite{GillSoul2}: 
$$\sum_{i=0}^{d}(-1)^{i}\length_R\big(\Tor_i^R\big(R/\phi^n(\mathfrak{m})R,S\big)\big)=0.$$
Since $f$ is a surjective morphism of local algebraic dynamical systems, we have 
$$f(\phi^n(\mathfrak{m})R)S=\psi^n(f(\mathfrak{m})S)S=\psi^n(\mathfrak{n})S.$$ Hence, there are $R$-module isomorphisms
$$\Tor^R_0\big(R/\phi^n(\mathfrak{m})R,S\big)\cong\left(R/\phi^n(\mathfrak{m})R\right)\otimes_RS\cong S/\psi^n(\mathfrak{n})S.$$ We then obtain 
\begin{equation}
\label{Equation:4.2}
\length_S\big(S/\psi^n(\mathfrak{n})S\big)=\sum_{i=1}^{d}(-1)^{i-1}\length_R\big(\Tor_i^R\big(R/\phi^n(\mathfrak{m})R,S\big)\big).
\end{equation}
The result follows by applying logarithm to both sides of (\ref{Equation:4.2}) and letting $n\rightarrow\infty$.
\end{proof}
We should note that the alternating sum appearing on the right-hand sides of~(\ref{Equation:4.1}) and~(\ref{Equation:4.2}) is the partial Euler characteristic $\chi_1^R\left(R/\phi^n(\mathfrak{m})R,S\right)$ with the notation of~\cite{Licht}).\par Theorem~\ref{Theorem:4.1} can be applied to any local algebraic dynamical system, in which the local ring is of equal characteristic, as described in the next example.
\begin{example}
\label{Example:4.2}
Let \((S,\mathfrak{n},\psi)\) be a local algebraic dynamical system and assume that $S$ is of equal characteristic and not regular. Suppose $\mathfrak{n}$ can be generated by $d$ elements. Let $\hat{S}$ be the $\mathfrak{n}$-adic completion of $S$ and let $\hat{\psi}\colon\hat{S}\rightarrow\hat{S}$ be the endomorphism induced by $\psi$. Then by Cohen's Structure Theorem there exists a surjective homomorphism $\pi\colon R=k\llbracket X_1,\ldots,X_d\rrbracket\twoheadrightarrow\hat{S}$, where $k$ is the residue field of $S$. By~\cite[Theorem~3]{MajMiaSzp} the endomorphism $\hat{\psi}$ can be lifted to an endomorphism of finite length $\phi\colon R\rightarrow R$ in such a way that $\pi\circ\phi=\hat{\psi}\circ\pi$. Since $S\rightarrow\hat{S}$ is flat, by~\cite[Proposition~20]{MajMiaSzp} we have $h_{\mathrm{loc}}(\psi)=h_{\mathrm{loc}}(\hat{\psi})$. Thus, letting $\mathfrak{m}$ be the maximal ideal of $R$, by Theorem~\ref{Theorem:4.1} the following equality holds:
$$
h_{\mathrm{loc}}(\psi)=\lim_{n\rightarrow\infty}\frac{1}{n}\log\Big(\sum_{i=1}^d(-1)^{i-1}\length_R\big(\Tor_i^R\big(R/\phi^n(\mathfrak{m})R,\hat{S}\big)\big)\Big).
$$
\end{example}
\section{Open problems}
\label{Section:5}
We list a couple of open problems here that are of particular interest to us.
\begin{problem}
\label{Problem:1}
In the context of Theorem~\ref{Theorem:3.3} (with or without assuming Cohen-Macaulayness of $S$), is
$h_t(\mathbb{L}\psi^{\varstar})=h_t(\mathbb{L}\phi^{\varstar})+h_t(\mathbb{L}{\overline{\psi}}^{\varstar})$?
\end{problem}
\begin{problem}
\label{Problem:2}
Does Theorem~\ref{Theorem:2.4} extend to Cohen-Macaulay rings?
\end{problem}

\bibliographystyle{plain}
\bibliography{Myreferences3}

\end{document}